\documentclass{amsart}
\usepackage{amsfonts}
\usepackage{amsmath}
\usepackage{amssymb}
\usepackage{amsthm}
\usepackage{graphicx}
\newtheorem{theorem}{Theorem}
\newtheorem{corollary}{Corollary}
\newtheorem{example}{Example}
\newtheorem{lemma}{Lemma}
\newtheorem{proposition}{Proposition}

\numberwithin{equation}{section}
\theoremstyle{definition}
\newtheorem{definition}{Definition}

\begin{document}
\title{Periodic Orbits of Billiards on an Equilateral Triangle}
\author{Andrew M. Baxter}
\address{Department of Mathematics, Rutgers University, 110 Frelinghuysen Rd, Piscataway NJ 08854}
\email{Andrew.Baxter@gmail.com}
\author{Ron Umble}
\address{Department of Mathematics, Millersville University of Pennsylvania, Millersville, PA 17551}
\email{ron.umble@millersville.edu}
\thanks{The results in this paper appeared in the first author's
undergraduate thesis supervised by the second author. }
\date{March 12, 2007}

\maketitle

\begin{abstract}
Using elementary methods, we find, classify and count the classes of
periodic orbits of a given period on an equilateral triangle. A periodic
orbit is either primitive or some iterate of a primitive orbit. Every
periodic orbit with odd period is some odd iterate of Fagnano's period 3.
Let $\mu$ denote the M\"{o}bius function. For each $n\in\mathbb{N}$, there
are exactly $\sum_{d\mid n}\mu(d)\left( \lfloor\frac{n/d+2}{2}\rfloor-\lfloor
\frac {n/d+2}{3}\rfloor\right) $ classes of primitive orbits with period $2n$.
\end{abstract}

\section{INTRODUCTION}

The trajectory of a billiard ball in motion on a frictionless billiards
table is completely determined by its initial position, direction, and speed.
When the ball strikes a bumper, we assume that the angle of incidence equals
the angle of reflection. Once released, the ball continues indefinitely
along its trajectory with constant speed unless it strikes a vertex, at
which point it stops. If the ball returns to its initial position with its
initial velocity direction, it retraces its trajectory and continues to do
so repeatedly; we call such trajectories \emph{periodic}. Nonperiodic
trajectories are either \emph{infinite} or \emph{singular}; in the later
case the trajectory terminates at a vertex. 
 
More precisely, think of a billiards table as a plane region $R$ bounded by a polygon $G$.    A \emph{nonsingular trajectory on} $G$ is a piecewise linear constant speed curve $\alpha : \mathbb{R} \rightarrow R$, where $\alpha(t)$ is the position of the ball at time $t$. An \emph{orbit} is the restriction of some nonsingular
trajectory to a closed interval; this is distinct from the notion of ``orbit'' in
discrete dynamical systems.

A nonsingular trajectory $\alpha$ is \emph{periodic} if $\alpha(a+t)=\alpha(b+t)$
for some $a < b$ and all $t\in \mathbb{R}$; its restriction to $[a,b]$ is a
\emph{periodic orbit}. A periodic orbit retraces the same path exactly $n\geq 1$ times. If $n=1$, the orbit is \emph{primitive}; otherwise it is an $n$-\emph{fold iterate. }
If $\alpha$ is primitive, $\alpha ^{n}$ denotes its $n$-fold iterate.  The \emph{period} of a periodic orbit is the number of times the ball strikes a bumper as it travels along its trajectory. If $\alpha $ is primitive of period $k$, then $\alpha ^{n}$ has period $kn$. 

In this article we give a complete solution to the following billiards
problem:\ \textit{Find, classify, and count the classes of periodic orbits of
a given period on an equilateral triangle.} While periodic orbits are known
to exist on all nonobtuse and certain classes of obtuse triangles \cite{H-H}
, \cite{M}, \cite{S}, \cite{V-G-S}, existence in general remains a
long-standing open problem. The first examples of periodic orbits were
discovered by Fagnano in 1745. Interestingly, his orbit of period 3 on an
acute triangle, known as the \textquotedblleft Fagnano
orbit,\textquotedblright\ was not found as the solution of a billiards
problem, but rather as the triangle of least perimeter inscribed in a given
acute triangle. This problem, known as \textquotedblleft Fagnano's
problem,\textquotedblright\ is solved by the orthic triangle, whose vertices
are the feet of the altitudes of the given triangle (see Figure 1). The
orthic triangle is a periodic trajectory since its angles are bisected by the
altitudes of the triangle in which it is inscribed; the proof given by
Coxeter and Greitzer \cite{C} uses exactly the \textquotedblleft
unfolding\textquotedblright\ technique we apply below. Coxeter credits this
technique to H. A. Schwarz and mentions that Frank and F. V. Morley \cite{M1}
extended Schwarz's treatment on triangles to odd-sided polygons. For a
discussion of some interesting properties of the Fagnano orbit on any acute
triangle, see \cite{G1}.

\begin{figure}[h]
\begin{center}
\includegraphics[width=3.3in]{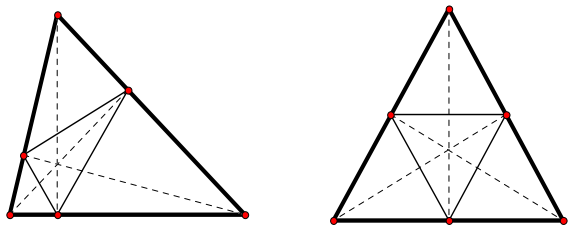}
\end{center}
\caption{Fagnano's period 3 orbit.}
\end{figure}

Much later, in 1986, Masur \cite{M} proved that every $\emph{rational}$
polygon (one whose interior angles are rational multiples of $\pi$) admits
infinitely many periodic orbits with distinct periods, but he neither
constructed nor classified them. A year later Katok \cite{K} proved that the
number of periodic orbits of a given period grows subexponentially.
Existence results on various polygons were compiled by Tabachnikov \cite{T}
in 1995.

This article is organized as follows: In Section 2 we introduce an
equivalence relation on the set of all periodic orbits on an equilateral
triangle and prove that every orbit with odd period is an odd iterate of
Fagnano's orbit. In Section 3 we use techniques from analytic geometry to
identify and classify all periodic orbits. The paper concludes with Section 4,
in which we derive two counting formulas: First, we establish a bijection
between classes of orbits with period $2n$ and partitions of $n$ with 2 or 3
as parts and use it to show that there are $\mathcal{O}(n)=\lfloor \frac{n+2
}{2}\rfloor -\lfloor \frac{n+2}{3}\rfloor $ classes of orbits with period $
2n $ (counting iterates). Second, we show that there are $\mathcal{P}
(n)=\sum_{d\mid n}\mu (d)\mathcal{O}\left( n/d\right) $ classes of primitive
orbits with period $2n,$ where $\mu $ denotes the M\"{o}bius function.

\section{ORBITS AND TESSELLATIONS}

Consider an equilateral triangle $\bigtriangleup ABC.$ We begin with some
key observations.

\begin{proposition}
\label{range}Every nonsingular trajectory strikes some side of $\bigtriangleup
ABC$ with an angle of incidence in the range $30^{\circ}\leq\theta
\leq60^{\circ}.$
\end{proposition}

\begin{proof}
Given a nonsingular trajectory $\alpha ,$ choose a point $P_{1}$ at which $%
\alpha $ strikes $\bigtriangleup ABC$ with angle of incidence $\theta _{1}.$
If $\theta _{1}$ lies in the desired range, set $\theta =\theta _{1}$.
Otherwise, let $\alpha _{1}$ be the segment of $\alpha $ that connects $P_{1}
$ to the next strike point $P_{2}$ and label the vertices of $\bigtriangleup
ABC$ so that $P_{1}$ is on side $\overline{AC}$ and $P_{2}$ is on side $%
\overline{BC}$ (see Figure 2). If $0^{\circ }<\theta _{1}<30^{\circ },$ then $\theta
_{2}=m\angle P_{1}P_{2}B=\theta _{1}+60^{\circ }$ so that $60^{\circ
}<\theta _{2}<90^{\circ }$. Let $\alpha _{2}$ be the segment of $\alpha $
that connects $P_{2}$ to the next strike point $P_{3}.$ Then the angle of
incidence at $P_{3}$ satisfies $30^{\circ }<\theta _{3}<60^{\circ };$ set $%
\theta =\theta _{3}.$ If $60^{\circ }<\theta _{1}\leq 90^{\circ }$ and $%
\theta _{1}$ is an interior angle of $\bigtriangleup P_{1}P_{2}C,$ then the
angle of incidence at $P_{2}$ is $\theta _{2}=m\angle P_{1}P_{2}C=120^{\circ
}-\theta _{1}$ and satisfies $30^{\circ }\leq \theta _{2}<60^{\circ }$; set 
$\theta =\theta _{2}.$ But if $60^{\circ }<\theta _{1}\leq
90^{\circ }$ and $\theta _{1}$ is an exterior angle of $\bigtriangleup
P_{1}P_{2}C,$ then the angle of incidence at $P_{2}$ is $\theta _{2}=m\angle
P_{1}P_{2}C=\theta _{1}-60^{\circ },$ in which case $0^{\circ }<\theta
_{2}\leq 30^{\circ }.$ If  $\theta _{2}=30^{\circ }$ set $\theta =\theta
_{2};$ otherwise continue as above until $30^{\circ }<\theta _{4}<60^{\circ }
$ and set $\theta =\theta _{4}.$
\end{proof}

\begin{figure}[h]
\begin{center}
\includegraphics[width=4.7 in]{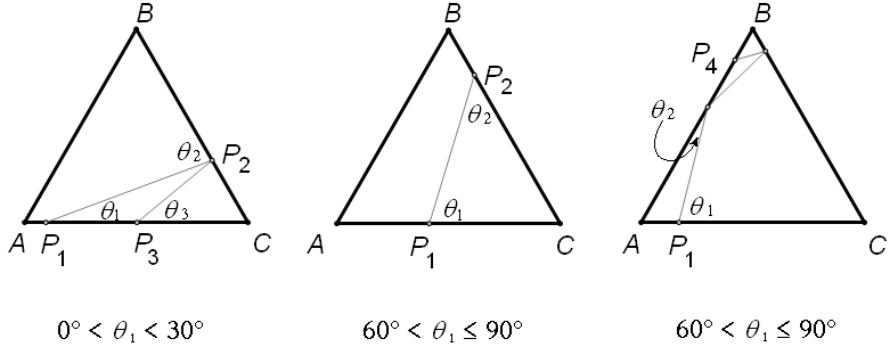}
\end{center}
\caption{Incidence angles in the range $30^{\circ}\leq\theta\leq 60^{\circ}.$}
\end{figure}

Let $\alpha $ be an orbit of period $n$ on $\bigtriangleup ABC$ oriented 
so that $\overline{BC}$ is horizontal. Since
Proposition \ref{range} applies equally well to periodic orbits, choose a
point $P$ at which $\alpha $ strikes $\bigtriangleup ABC$ with angle of
incidence in the range $30^{\circ }\leq \theta \leq 60^{\circ }.$ If
necessary, relabel the vertices of $\bigtriangleup ABC$, change initial
points, and reverse the parameter so that side $\overline{BC}$ contains $P,$ $%
\alpha $ begins and ends at $P$, and the components of $\alpha ^{\prime }$ as
the ball departs from $P$ are positive. Let $\mathcal{T}$ be a regular
tessellation of the plane by equilateral triangles, each congruent to $%
\bigtriangleup ABC$, and positioned so that one of its families of parallel
edges is horizontal. Embed $\bigtriangleup ABC$ in $\mathcal{T}$ so that its
base $\overline{BC}$ is collinear with a horizontal edge of $\mathcal{T}$.
Let $\alpha _{1},\alpha _{2},\ldots ,\alpha _{n}$ denote the directed
segments of $\alpha ,$ labelled sequentially; then $\alpha _{1}$ begins at $P
$ and terminates at $P_{1}$ on side $s_{1}$ of $\bigtriangleup ABC$ with
angle of incidence $\theta _{1}$. Let $\sigma _{1}$ be the reflection in the
edge of $\mathcal{T}$ containing $s_{1}.$ Then $\alpha _{1}$ and $\sigma
_{1}(\alpha _{2})$ are collinear segments and $\sigma _{1}(\alpha )$ is a
periodic orbit on $\sigma _{1}(\bigtriangleup ABC),$ which is the basic
triangle of $\mathcal{T}$ sharing side $s_{1}$ with $\bigtriangleup ABC$.
Follow $\sigma _{1}(\alpha _{2})$ from $P_{1}$ until it strikes side $s_{2}$
of $\sigma _{1}(\bigtriangleup ABC)$ at $P_{2}$ with incidence angle $\theta
_{2}$. Let $\sigma _{2}$ be the reflection in the edge of $\mathcal{T}$
containing $s_{2};$ then $\alpha _{1}$, $\sigma _{1}(\alpha _{2})$ and $%
\left( \sigma _{2}\sigma _{1}\right) (\alpha _{3})$ are collinear segments
and $\left( \sigma _{2}\sigma _{1}\right) (\alpha )$ is a periodic orbit on $%
\left( \sigma _{2}\sigma _{1}\right) (\bigtriangleup ABC)$. Continuing in
this manner for $n-1$ steps, let $\theta _{n}$ be the angle of incidence at $%
Q=\left( \sigma _{n-1}\sigma _{n-2}\cdots \sigma _{1}\right) \left( P\right)
.$ Then $\theta _{1},\theta _{2},\ldots ,\theta _{n}$ is a sequence of
incidence angles with $30^{\circ }\leq \theta _{n}\leq 60^{\circ },$ and $%
\alpha _{1},\sigma _{1}(\alpha _{2}),\ldots ,\left( \sigma _{n-1}\sigma
_{n-2}\cdots \sigma _{1}\right) \left( \alpha _{n}\right) $ is a sequence of
collinear segments whose union is the directed segment from $P$ to $Q.$
Using the notation in \cite{Ma}, let $\underline{PQ}$ denote the directed
segment from $P$ to $Q.$ Then $\underline{PQ}$ has the same length as $%
\alpha $ and enters and exits the triangle $\left( \sigma _{i}\cdots \sigma
_{1}\right) \left( \bigtriangleup ABC\right) $ with angles of incidence $%
\theta _{i}$ and $\theta _{i+1}.$ We refer to $\underline{PQ}$ as an \emph{\
unfolding }of $\alpha $ and to $\theta _{n}$ as its \emph{representation
angle}.

\begin{proposition}
\label{incidence}A periodic orbit strikes the sides of $\bigtriangleup ABC$
with at most three incidence angles, exactly one of which lies in the range $
30^{\circ}\leq\theta\leq60^{\circ}.$ In fact, exactly one of the following
holds:

\begin{enumerate}
\item All incidence angles measure $60^{\circ}.$

\item There are exactly two distinct incidence angles measuring $30^{\circ}$
and $90^{\circ}$.

\item There are exactly three distinct incidence angles $\phi$, $\theta$, and 
$\psi$ such that $0^{\circ}<\phi<30^{\circ}<\theta<60^{\circ}<\psi<90^{\circ
}.$
\end{enumerate}
\end{proposition}

\begin{proof}
Let $\alpha $ be a periodic orbit and let $\underline{PQ}$ be an unfolding.
By construction, $\underline{PQ}$ cuts each horizontal edge of $\mathcal{T
}$ with angle of incidence in the range $30^{\circ }\leq \theta \leq
60^{\circ }.$ Consequently, $\underline{PQ}$ cuts a left-leaning edge of $
\mathcal{T}$ with angle of incidence $\phi =120^{\circ }-\theta $ and
cuts a right-leaning edge of $\mathcal{T}$ with angle of incidence $\psi
=60^{\circ }-\theta $ (see Figure 3). In particular, if $\theta =60^{\circ
}, $ $\underline{PQ}$ cuts only left-leaning and horizontal edges, and
all incidence angles are equal. In this case, $\alpha $ is either the
Fagnano orbit, a primitive orbit of period 6, or some iterate of these. If $
\theta =30^{\circ },$ then $\phi =90^{\circ }$ and $\psi =30^{\circ },$ and $
\alpha $ is either primitive of period $4$ or some iterate thereof (see
Figure 4). When $30^{\circ }<\theta <60^{\circ },$ clearly $0^{\circ
}<\phi <30^{\circ }$ and $60^{\circ }<\psi <90^{\circ }.$
\end{proof}

\begin{corollary}
\label{parallel}Any two unfoldings of a periodic orbit are parallel.
\end{corollary}

\begin{figure}[h]
\begin{center}
\includegraphics[height=1.1 in]{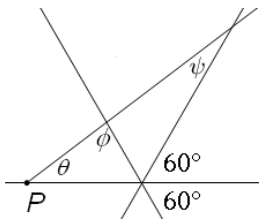}
\end{center}
\caption{Incidence angles $\protect\theta,$ $\protect\phi$, and $\protect
\psi.$}
\end{figure}

\begin{figure}[h]
\begin{center}
\includegraphics[width=5.3 in]{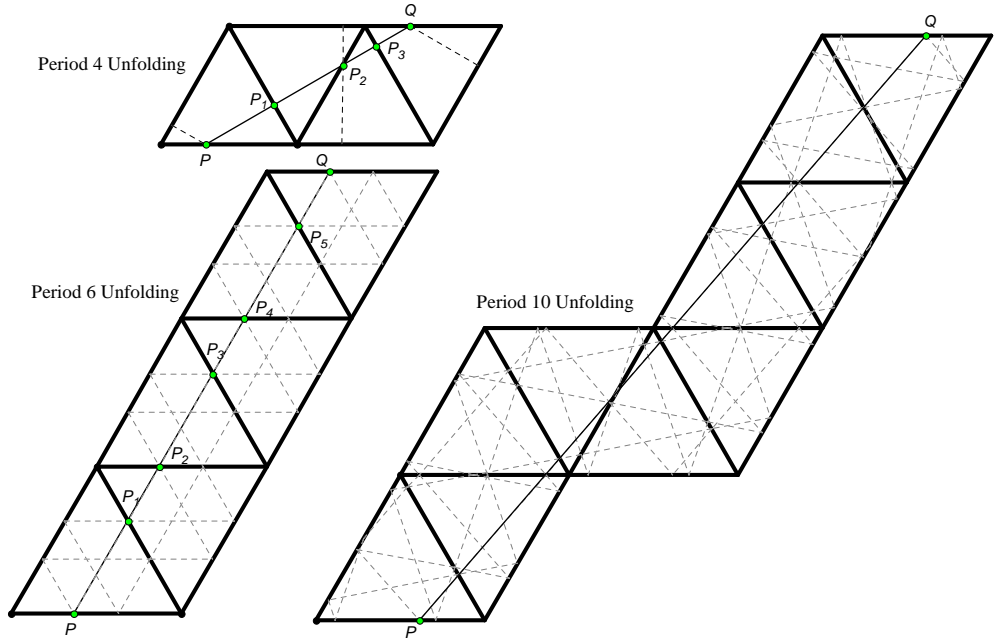}
\end{center}
\caption{Unfolded orbits of period 4, 6, and 10.}
\end{figure}

Our next result plays a pivotal role in the classification of orbits.

\begin{theorem}
\label{period}If an unfolding of a periodic orbit $\alpha $ terminates on a
horizontal edge of $\mathcal{T}$, then $\alpha $ has even period.
\end{theorem}

\begin{proof}
Let $\underline{PQ}$ be an unfolding of $\alpha .$ Then both $P$ and $Q$ lie
on horizontal edges of $\mathcal{T}$, and the basic triangles of $\mathcal{T}
$ cut by $\underline{PQ}$ pair off and form a polygon of rhombic tiles
containing $\underline{PQ}$ (see Figure 5). As the path $\underline{PQ}$
traverses this polygon, it enters each rhombic tile through an edge, cuts
a diagonal of that tile (collinear with a left-leaning edge of $\mathcal{T}$), and exits through another edge. Since each exit edge of one tile is the
entrance edge of the next and the edge containing $P$ is identified with the
edge containing $Q$, the number of distinct edges of $\mathcal{T}$ cut by $
\underline{PQ}$ is twice the number of rhombic tiles. It follows that $\alpha $ has
even period.
\end{proof}

\begin{figure}[h]
\begin{center}
\includegraphics[width=3.4 in]{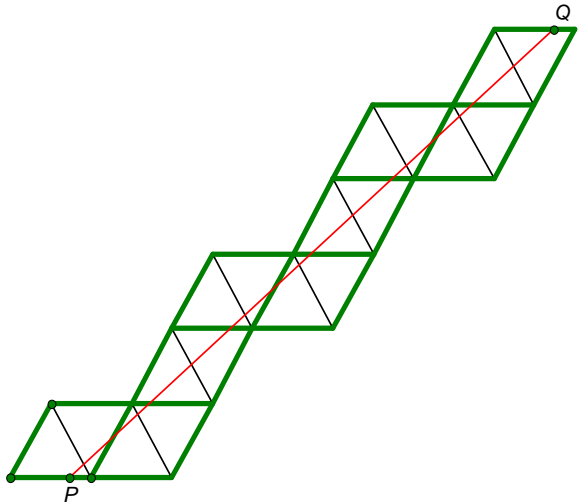}
\end{center}
\caption{A typical rhombic tiling.}
\end{figure}

Let $\gamma $ denote the Fagnano orbit.

\begin{theorem}
\label{horiz}If $\alpha $ is a periodic orbit and $\alpha \neq \gamma
^{2k-1} $ for all $k\geq 1,$ then every unfolding of $\alpha $ terminates
on a horizontal edge of $\mathcal{T}$.
\end{theorem}

\begin{proof}
We prove the contrapositive. Suppose there is an unfolding $\underline{PQ}$
of $\alpha $ that does not terminate on a horizontal edge of $\mathcal{T}$.
Let $\theta $ be the angle of incidence at $Q;$ then $\theta $ is also the
angle of incidence at $P$ and $\theta \in \left\{ 30^{\circ },60^{\circ
}\right\} $ by the proof of Proposition \ref{incidence}. But if $\theta
=30^{\circ }$, then $\alpha $ is some iterate of the period 4 orbit whose
unfoldings terminate on a horizontal edge of $\mathcal{T}$ (see Figure 4). So $\theta =60^{\circ }.$ But $\alpha $ is neither an iterate of a period 6
orbit nor an even iterate of $\gamma $ since their unfoldings also terminate
on a horizontal edge of $\mathcal{T}$ (see Figure 4). It follows that $
\alpha =\gamma ^{2k-1}$ for some $k\geq 1$.
\end{proof}

Combining the contrapositives of Theorems \ref{period} and \ref{horiz} we
obtain the following characterization:

\begin{corollary}
\label{odd}If $\alpha $ is an orbit with odd period, then $\alpha =\gamma
^{2k-1}$ for some $k\geq 1,$ in which case the period is $6k-3.$
\end{corollary}

Let $\alpha $ be an orbit with even period and let $\underline{PQ}$ be an
unfolding. Let $G$ be the group generated by all reflections in the edges of 
$\mathcal{T}$. Since the action of $G$ on $\overline{BC}$ generates a
regular tessellation $\mathcal{H}$ of the plane by hexagons, $\alpha $
terminates on some horizontal edge of $\mathcal{H}$. As in the definition of an unfolding, let $\sigma_1,\sigma_2,\ldots,\sigma_{n-1}$ be the reflections in the lines of $\mathcal{T}$ cut by $\underline{PQ}$ (in order) and $\sigma_n$ be the reflection in the line of $\mathcal{T}$ containing $Q$.  Then the composition $f=\sigma_{n}\sigma _{n-1}\cdots \sigma _{1}$ maps $P$ to $Q$ and maps the hexagon whose base $\overline{BC}$ contains $P$ to the hexagon whose base $\overline{B'C'}$ contains $Q$.  Then $n$ (the period of $\alpha$) is even and $f$ is either a translation by vector $\overrightarrow{PQ}$ or a rotation of $120^{\circ }$ or $240^{\circ }$. But $\overline{BC}\Vert \overline{B'C'}$ so $f$ is a translation and the position of $Q$ on $\overline{B'C'}$ is exactly the same as the position of $P$ on $
\overline{BC}$.

Periodic orbits represented by horizontal translations of
an unfolding $\underline{PQ}$ are generically distinct, but have the same
length and incidence angles (up to permutation) as $\alpha $. Hence it is
natural to think of them as equivalent.

\begin{definition}
Periodic orbits $\alpha $ and $\beta $ are {\it equivalent} if there exist respective unfoldings $\underline{PQ}$ and $\underline{RS}$ and a
horizontal translation $\tau $ such that $\underline{RS}=\tau \left( 
\underline{PQ}\right) $. The symbol $\left[ \alpha \right] $ denotes the
equivalence class of $\alpha .$ The {\it period} of a class $\left[
\alpha \right] $ is the period of its elements; a class is {\it even}
if and only if it has even period.
\end{definition}

\begin{figure}[h]
\begin{center}
\includegraphics[width=2.7 in]{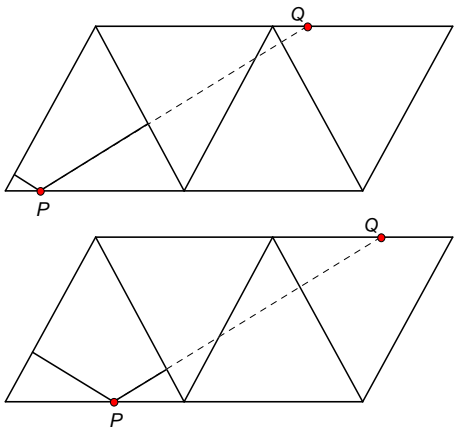}
\end{center}
\caption{Unfoldings of equivalent period 4 orbits.}
\end{figure}

Consider an unfolding $\underline{PQ}$ of a periodic orbit $\alpha$. If $
\left[ \alpha \right] $ is even, let $R$ be a point on $\overline{BC}$ and
let $\tau $ is the translation from $P$ to $R$. We say that the point $R$ is 
\emph{singular for $[\alpha ]$} if $\tau (\underline{PQ})$ contains a vertex
of $\mathcal{T}$; then $\tau (\underline{PQ})$ is an unfolding of a periodic
orbit whenever $R$ is non-singular for $[\alpha]$.  Furthermore, $\alpha  $ strikes $\underline{BC}$ at finitely many points and at most finitely many points on $\overline{BC}$ are singular for $\left[ \alpha \right]$.  Therefore $[\alpha]$ has cardinality $\mathfrak{c}$ (the cardinality of an interval). On the other hand, Corollary \ref{odd} tells us that an orbit of odd period is $\gamma^{2k-1}$ for some $k\geq 1$.  But if $k \neq \ell$, then $\gamma^{2k-1}$ and $\gamma^{2\ell-1}$ have different periods and cannot be equivalent. Therefore $[\gamma^{2k-1}]$ is a singleton class for each $k$. We have proved:

\begin{proposition}
\label{parity}The cardinality of a class is determined by its parity; in
fact, $\alpha $ has odd period if and only if $\left[ \alpha \right] $ is a singleton class.
\end{proposition}

Proposition \ref{parity} and Corollary \ref{odd} completely classify orbits with odd period. The remainder of this article considers orbits with even period. Our strategy is to represent the classes of all such orbits as lattice points in some ``fundamental region," which we now define. First note that any two unfoldings whose terminal points lie on the same horizontal edge of $\mathcal{H}$ are equivalent. Since $\mathcal{H}$ has countably many horizontal edges, there are countably many even classes of orbits.  Furthermore, since at most finitely many points in $\overline{BC}$ are singular for each even class, there is a point $O$ on $\overline{BC}$ other than the midpoint that is nonsingular for every class.  Therefore, given an even class $[\alpha]$, there is a point $S$ and an element $x \in \left[\alpha\right]$ such that $\underline{OS}$ is an unfolding of $x$.  Note that if $\underline{PQ}$ is an unfolding of $\alpha$, then $\underline{OS}$ is the horizontal translation of $\underline{PQ}$ by $\overrightarrow{PO}$.  Therefore $\alpha$ uniquely determines the point $S$, denoted henceforth by $S_{\alpha}$, and we refer to $\underline{OS_{\alpha}}$ as the \emph{fundamental unfolding of $[\alpha ]$}. The \emph{fundamental region at O}, denoted by $\Gamma _{O}$, is the polar region $30^{\circ }\leq \theta \leq 60^{\circ }$ centered at $O$; the points $S_{\alpha}$ are called \emph{lattice points} of $\Gamma _{O}.$

Since $O$ is not the midpoint of $\overline{BC}$, odd iterates of
Fagnano's orbit $\gamma$ have no fundamental unfoldings. On the other hand,
the fundamental unfolding of $\gamma ^{2n}$ represents the $n$-fold iterate of
a primitive period 6 orbit. Nevertheless, with the notable exception of $[ \gamma ^{2}] $, ``primitivity" is a property common to all orbits of the
same class (see Figure 7). Indeed, the fundamental unfolding of $[\gamma ^{2}]$ represents a primitive orbit. So we define a \emph{primitive
class} to be either $\left[ \gamma ^{2}\right] $ or a class of primitives.

\begin{figure}[h]
\begin{center}
\includegraphics[width=2.5 in]{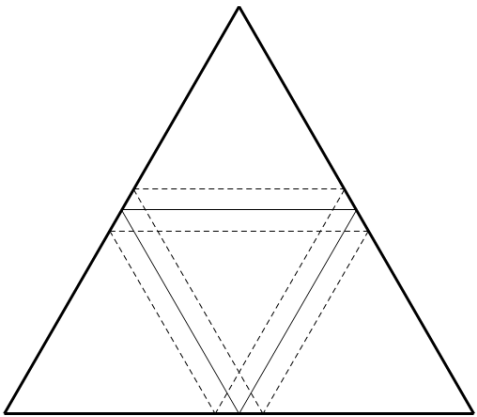}
\end{center}
\caption{The Fagnano orbit and an equivalent period 6 orbit (dotted).}
\end{figure}

To complete the classification, we must determine exactly which directed
segments in $\Gamma _{O}$ with initial point $O$ represent orbits with even
period. We address this question in the next section.

\section{ORBITS AND RHOMBIC COORDINATES}

In this section we introduce the analytical structure we need to complete
the classification and to count the distinct classes of orbits of a given
even period. Expressing a fundamental unfolding $\underline{OS}$ as a vector 
$\overrightarrow{OS}$ allows us to exploit the natural rhombic coordinate
system given by $\mathcal{T}$. Let $O$ be the origin and take the $x$-axis
to be the horizontal line containing it. Take the $y$-axis to be the line
through $O$ with inclination $60^{\circ }$ and let $BC$ be the unit of
length (see Figure 8). Then in rhombic coordinates
\begin{equation*}
\Gamma _{O}=\left\{ \left( x,y\right) \text{ }|\text{ }0\leq x\leq y\right\}.
\end{equation*}

\begin{figure}[h]
\begin{center}
\includegraphics[width=2.9 in]{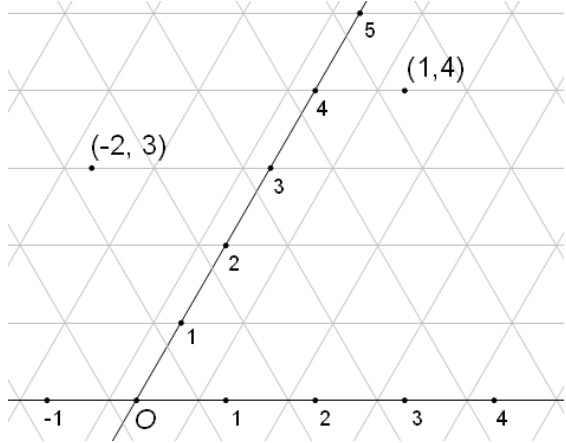}
\end{center}
\caption{Rhombic coordinates.}
\end{figure}

Since the period of $\left[ \alpha \right] $ is twice the number of rhombic
tiles cut by $\underline{OS_{\alpha}},$ and the rhombic
coordinates of $S_{\alpha}$ count these rhombic tiles, we can strengthen Theorem \ref
{period}:

\begin{corollary}
\label{corperiod}If $S_{\alpha}=\left( x,y\right) $, then $\alpha$ has period $2(x+y).$
\end{corollary}

Points in the integer sublattice $\mathcal{L}$ of points on the horizontals of $\mathcal{H}$ that are images of $O$ under the action of $G$ have the following simple characterization: Let $H$ be the hexagon of $\mathcal{H}$ with base $\overline{BC}$, and let $\tau_1$ and $\tau_2$ denote the translations by the vectors $\left( 1,1 \right)$ and $\left( 0,3 \right)$, respectively.  Then the six hexagons adjacent to $H$ are its images $\tau_{2}^{b}\tau _{1}^{a}\left( H\right) $, $(a,b) \in \left\{ \pm (1,0) ,\pm (1,-1) ,\pm (2,-1) \right\}$. Inductively, if $H'$ is any hexagon of $\mathcal{H}$, then $H'=\tau _{2}^{b}\tau _{1}^{a}\left( H\right) $ for some $a,b\in \mathbb{Z}$. Note that $a(1,1)+b(0,3)$ defines the translation $\tau_{2}^{b}\tau _{1}^{a}$. Hence $\mathcal{L}$ is generated by the vectors $(1,1)$ and $(0,3)$ and it follows that $(x,y)\in \mathcal{L}$ if and only if $x\equiv y\pmod{3}$.

Now recall that if $\underline{PQ}$ is an unfolding, then $Q$ lies on a
horizontal of $\mathcal{H}$. Hence $\underline{OS}$ is a fundamental unfolding if and
only if $S\in \mathcal{L}\cap \Gamma _{O}-O$ if and only if $S \in \{(x,y) \in
\mathbb{Z}^{2}\cap \Gamma _{O} \mid x\equiv y(\text{\textrm{mod }}3),x+y=n\}$. We have proved:

\begin{theorem}
\label{correspond} Given an even class $[\alpha]$, let $(x,y)_{\alpha}=S_{\alpha}$. There is a bijection 
\begin{equation*}
\Phi:\{[\alpha ] \mid [\alpha ]\text{ has period }2n\}\rightarrow \{(x,y)\in
\mathbb{Z}^{2}\cap \Gamma _{O} \mid x\equiv y (\text{\textrm{mod }}3),x+y=n\}
\end{equation*}
given by $\Phi\left([\alpha]\right)=(x,y)_{\alpha}$.
\end{theorem}

Taken together, Proposition \ref{parity}, Corollary \ref{corperiod}
and Theorem \ref{correspond} classify all periodic orbits on an equilateral
triangle.

\begin{theorem}
\label{classify}(Classification) Let $\alpha $ be a periodic orbit on an
equilateral triangle.

\begin{enumerate}
\item If $\alpha $ has period $2n$, then $\left[ \alpha \right] $ has
cardinality $\mathfrak{c}$ and contains exactly one representative whose unfolding $
\underline{OS}$ satisfies $S=\left( x,y\right)$, $0\leq x\leq y,$ 
$x\equiv y\pmod{3}$, and $x+y=n$.

\item Otherwise, $\alpha =\gamma ^{2k-1}$ for some $k\geq 1,$ in which case
its period is $6k-3.$ 
\end{enumerate}
\end{theorem}

In view of Theorem \ref{correspond}, we may count classes of orbits of a
given period $2n$ by counting integer pairs $(x,y)$ such that $0\leq x\leq
y, $ $x\equiv y\pmod{3}$ and $x+y=n.$ This is the objective of the next and
concluding section.
\begin{figure}[h]
\begin{center}
\includegraphics[width=5 in]{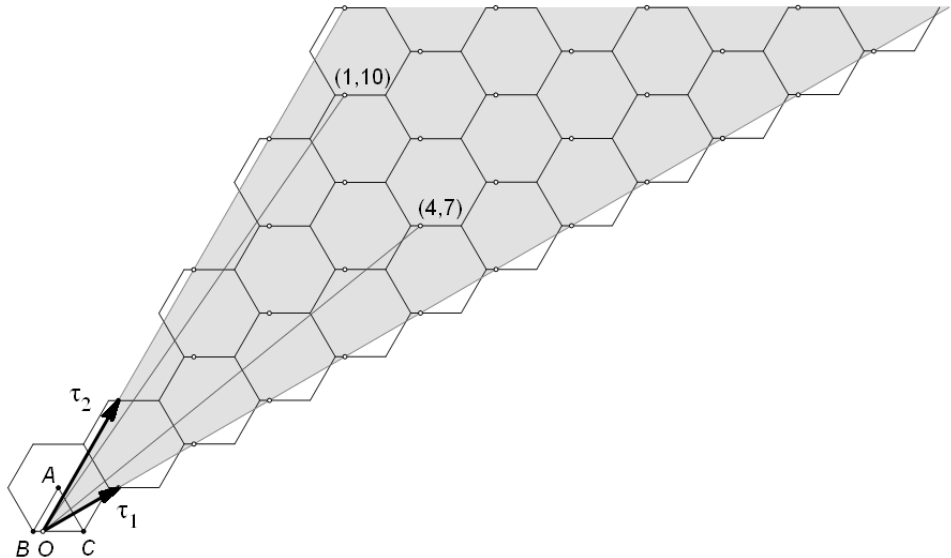}
\end{center}
\caption{Translated images of $O$ in $\Gamma _{O}$ and unfoldings of period
22 orbits.}
\end{figure}

\section{ORBITS AND INTEGER PARTITIONS}

We will often refer to an ordered pair $\left( x,y\right) $ as an
``orbit" when we mean the even class of orbits to which it corresponds. Two questions arise: (1) Is there an orbit
with period $2n$ for each $n\in \mathbb{N}$? (2) If so, exactly how many distinct classes of orbits with period $2n$ are there?

\begin{figure}[h]
\begin{center}
\includegraphics[width=4.9 in]{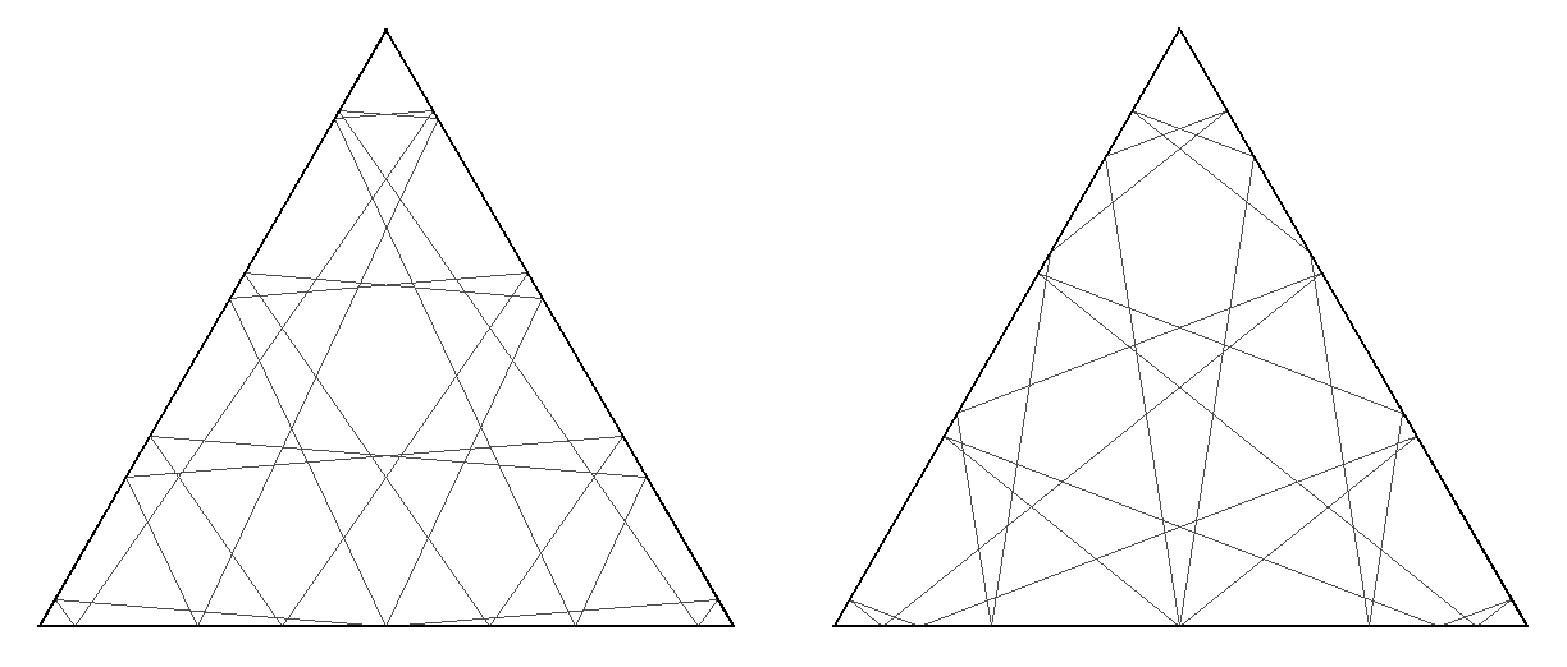}
\end{center}
\caption{Period 22 orbits $(1,10)$ (left) and $(4,7)$ (right).}
\end{figure}

If we admit iterates, question (1) has an easy answer.  Clearly there are no period 2 orbits since no two sides of $\bigtriangleup ABC$ are parallel --- alternatively, if $(a,b)$ is a solution of the system $x\equiv y \pmod{3}$ and $x+y=1$, either $a$ or $b$ is negative.  For each $n>1,$ the orbit
\begin{equation*}
\alpha =\left\{ 
\begin{array}{ll}
(\frac{n}{2},\frac{n}{2}), & n\text{ even} \\ 
&  \\ 
(\frac{n-1}{2}-1,\frac{n-1}{2}+2), & n\text{ odd.}
\end{array}
\right.
\end{equation*}

\noindent has period $2n.$ Note that the period 22 orbits $(1,10)$ and $(4,7)$ are not
equivalent since they have different lengths and representation angles (see
Figures 9 and 10).

To answer to question (2), we reduce the problem to counting partitions by
constructing a bijection between classes of orbits with period $2n$ and
partitions of $n$ with 2 and 3 as parts. For a positive integer $n$, a \emph{
partition} of $n$ is a nonincreasing sequence of nonnegative integers whose
terms sum to $n$. Such a sequence has finitely many nonzero terms, called
the \emph{parts}$,$ followed by infinitely many zeros. Thus, we seek pairs
of nonnegative integers $(a,b)$ such that $n=2a+3b$. The reader can easily
prove:

\begin{lemma}
\label{biject}For each $n\in \mathbb{N},$ let 
\begin{equation*}
X_{n}=\left\{ \left( x,y\right) \in \mathbb{Z}^{2} \mid 0\leq
x\leq y,\text{ } x\equiv y\hspace*{-0.1in}\pmod{3},\text{ } x+y=n\right\} \text{and}
\end{equation*}
\begin{equation*}
Y_{n}=\{(a,b)\in \mathbb{Z}^{2} \mid a,b\geq 0\text{ and }
2a+3b=n\}.
\end{equation*}
The function $\varphi :Y_{n}\rightarrow X_{n}$ given by $\varphi \left(
a,b\right) =\left( a,a+3b\right) $ is a bijection.
\end{lemma}

Combining Theorem \ref{classify} and Lemma \ref{biject}, we have:

\begin{corollary}
For each $n\in \mathbb{N}$, there is a bijection between period $2n$ orbits and the partitions of $n$ with 2 and 3 as parts.
\end{corollary}

Counting partitions of $n$ with specified parts is well understood (e.g.,
Sloane's A103221, \cite{OEIS}). The number of partitions of $n$ with $2$ and 
$3$ as parts is the coefficient of $x^{n}$ in the generating function 
\begin{eqnarray*}
f(x)&=&\sum_{n=0}^{\infty} \mathcal{O}(n)x^n \\
&=&(1+x^{2}+x^{4}+x^{6}+\cdots )(1+x^{3}+x^{6}+x^{9}+\cdots )\\
&=&\frac{1}{(1-x^{2})(1-x^{3})}.
\end{eqnarray*}
To compute this coefficient, let $\omega $ be a primitive cube root of unity
and perform a partial fractions decomposition. Then 
\begin{eqnarray*}
f(x)&=&\frac{1}{4(1+x)}+\frac{1}{4(1-x)}+\frac{1}{6(1-x)^{2}}+\frac{1}{9}
\left( \frac{1+2\omega }{\omega -x}+\frac{1+2\omega ^{2}}{\omega ^{2}-x}
\right) \\ 
&=&\frac{1}{4}\sum_{n=0}^{\infty }(-1)^{n}x^{n}+\frac{1}{4}\sum_{n=0}^{\infty
}x^{n}+\frac{1}{6}\sum_{n=0}^{\infty }(n+1)x^{n}\\
&& +\frac{1}{9} \sum_{n=0}^{\infty }(\omega ^{2n+2}+2\omega ^{2n}+\omega ^{n+1}+2\omega^{n})x^{n},
\end{eqnarray*}
and we have 
\[
\mathcal{O}(n) = \frac{(-1)^{n}}{4}+\frac{n}{6}+\frac{5}{12}+\frac{1}{9}
\left( \omega ^{2n+2}+2\omega ^{2n}+\omega ^{n+1}+2\omega ^{n}\right).
\]

By easy induction arguments, one can obtain the following simpler formulations (see \cite{OEIS}):

\begin{theorem}
\label{OCount}The number of distinct classes of period $2n$ is exactly 
\begin{eqnarray*}
\mathcal{O}(n) &=&\left\{ 
\begin{array}{ll}
\lfloor \frac{n}{6}\rfloor , & n\equiv 1 \pmod{6} \\ 
\lfloor \frac{n}{6}\rfloor +1, & \text{\rm otherwise}
\end{array}
\right. \\
&&\\
&=&\left\lfloor \frac{n+2}{2}\right\rfloor -\left\lfloor \frac{
n+2}{3}\right\rfloor .
\end{eqnarray*}
\end{theorem}

Let us refine this counting formula by counting only primitives. For every
divisor $d$ of $n$, the $(n/d)$-fold iterate of a primitive period $2d$
orbit has period $2n$. Hence, if $\mathcal{P}(n)$ denotes the number of
primitive classes of period $2n$, then
\begin{equation*}
\mathcal{O}(n)=\sum_{d|n}\mathcal{P}(d).
\end{equation*}
A formula for $\mathcal{P}(n)$ is a direct consequence of the M\"{o}bius
inversion formula (see \cite{R}). The M\"{o}bius function $\mu :\mathbb{
N\rightarrow }\left\{ -1,0,1\right\} $ is defined by
\begin{equation*}
\mu \left( d\right) =\left\{ 
\begin{array}{ll}
1, & d=1 \\ 
\left( -1\right) ^{r}, & d=p_{1}p_{2}\cdots p_{r}\text{ for distinct primes }
p_{i} \\ 
0, & \text{otherwise.}
\end{array}
\right. 
\end{equation*}

\begin{theorem}
\label{primcount}
For each $n\in \mathbb{N}$, there are exactly 
\begin{equation*}
\mathcal{P}(n)=\sum_{d|n}\mu (d)\mathcal{O}(n/d)
\end{equation*}
primitive classes of period $2n$.
\end{theorem}

Theorems \ref{OCount} and \ref{primcount}, together with Example \ref{ex} below, imply:

\begin{corollary}
\label{SpecCases}$\mathcal{O}(n)=0$ if and only if $n=1$; $\mathcal{P}(n)=0$
if and only if $n=1,4,6,10.$
\end{corollary}

\begin{corollary}
The following are equivalent:

\begin{enumerate}
\item The integer $n$ is 1 or prime.

\item $\mathcal{P}(n)=\mathcal{O}(n).$

\item All classes of period $2n$ are primitive.
\end{enumerate}
\end{corollary}

\noindent Table 1 displays some values of $\mathcal{O}$ and $
\mathcal{P}$. The values $\mathcal{O}(4)=1$, $\mathcal{P}(4)=0,$ and $\mathcal{P}(2)=1$, for
example, indicate that the single class of period $8$ contains only $2$-fold
iterates of the primitive orbits in the single class of period $4$.

\begin{table}[ht]
\begin{center}
\begin{tabular}{|c|c|c|c|}
\hline
$n$ & $2n$ & $\mathcal{O}(n)$ & $\mathcal{P}(n)$ \\ \hline
1 & 2 & 0 & 0 \\ 
2 & 4 & 1 & 1 \\ 
3 & 6 & 1 & 1 \\ 
4 & 8 & 1 & 0 \\ 
5 & 10 & 1 & 1 \\ 
6 & 12 & 2 & 0 \\ 
7 & 14 & 1 & 1 \\ 
8 & 16 & 2 & 1 \\ 
9 & 18 & 2 & 1 \\ 
10 & 20 & 2 & 0 \\ 
11 & 22 & 2 & 2 \\ 
12 & 24 & 3 & 1 \\ 
13 & 26 & 2 & 2 \\ 
14 & 28 & 3 & 1 \\ 
15 & 30 & 3 & 1 \\ 
16 & 32 & 3 & 1 \\ 
17 & 34 & 3 & 3 \\ 
18 & 36 & 4 & 1 \\ 
19 & 38 & 3 & 3 \\ 
20 & 40 & 4 & 2 \\ 
21 & 42 & 4 & 2 \\ 
22 & 44 & 4 & 1 \\ 
23 & 46 & 4 & 4 \\ 
24 & 48 & 5 & 1 \\ 
25 & 50 & 4 & 3 \\ 
26 & 52 & 5 & 2 \\ 
27 & 54 & 5 & 3 \\ 
28 & 56 & 5 & 2 \\ 
29 & 58 & 5 & 5 \\ 
30 & 60 & 6 & 2 \\ \hline
\end{tabular}
\begin{tabular}{|c|c|c|c|}
\hline
$n$ & $2n$ & $\mathcal{O}(n)$ & $\mathcal{P}(n)$ \\ \hline
31 & 62 & 5 & 5 \\ 
32 & 64 & 6 & 3 \\ 
33 & 66 & 6 & 3 \\ 
34 & 68 & 6 & 2 \\ 
35 & 70 & 6 & 4 \\ 
36 & 72 & 7 & 2 \\ 
37 & 74 & 6 & 6 \\ 
38 & 76 & 7 & 3 \\ 
39 & 78 & 7 & 4 \\ 
40 & 80 & 7 & 2 \\ 
41 & 82 & 7 & 7 \\ 
42 & 84 & 8 & 2 \\ 
43 & 86 & 7 & 7 \\ 
44 & 88 & 8 & 4 \\ 
45 & 90 & 8 & 4 \\ 
46 & 92 & 8 & 3 \\ 
47 & 94 & 8 & 8 \\ 
48 & 96 & 9 & 3 \\ 
49 & 98 & 8 & 7 \\ 
50 & 100 & 9 & 4 \\ 
51 & 102 & 9 & 5 \\ 
52 & 104 & 9 & 4 \\ 
53 & 106 & 9 & 9 \\ 
54 & 108 & 10 & 3 \\ 
55 & 110 & 9 & 6 \\ 
56 & 112 & 10 & 4 \\ 
57 & 114 & 10 & 6 \\ 
58 & 116 & 10 & 4 \\ 
59 & 118 & 10 & 10 \\ 
60 & 120 & 11 & 2 \\ \hline
\end{tabular}
\end{center}
\caption{Sample Values for $\mathcal{O}(n)$ and $\mathcal{P}(n)$.}
\end{table}

We conclude with an example of a primitive class of period $2n$ for each $
n\in \mathbb{N}-\{1,4,6,10\}$. But first we need the following self-evident lemma:

\begin{lemma}
\label{Primitive}Given an orbit $(x,y)\in \Gamma _{O},$ let $d\in \mathbb{N}$
be the largest value such that $x/d\equiv y/d\pmod{3}$. Then $(x,y)$ is
primitive if and only if $d=1$; otherwise $(x,y)$ is a $d$-fold iterate of
the primitive orbit $(x/d,y/d)$.
\end{lemma}

Although $d$ is difficult to compute, it is remarkably easy to check for
primitivity.

\begin{theorem}
\label{PrimeCheck}An orbit $(x,y)\in \Gamma _{O}$ is primitive if and only
if either

\begin{enumerate}
\item $\gcd\left( x,y\right) =1$ or

\item $\left( x,y\right) =(3a,3b),$ $\gcd\left( a,b\right) =1$, and $a\not
\equiv b \pmod{3}$ for some $a,b\in \mathbb{N\cup}
\left\{ 0\right\} .$
\end{enumerate}
\end{theorem}

\begin{proof}
If $\gcd\left( x,y\right) =1$, the orbit $\left( x,y\right) $ is primitive.
On the other hand, if $(x,y)=(3a,3b),a\not \equiv b \pmod{3} $, and $\gcd\left( a,b\right) =1$ for some $a,b$, let $d$ be as in
Lemma \ref{Primitive}. Then $d\neq3$ since $a\not \equiv b \pmod{3}$. But $\gcd\left( a,b\right) =1$ implies $d=1$, so $(x,y)$ is also
primitive when (2) holds.

Conversely, given a primitive orbit $(x,y),$ let $c=\gcd(x,y).$ Then $cm=x\leq y=cn$ for some $m,n\in\mathbb{N\cup}\{0\};$ thus $m\leq n,$ $\gcd(m,n)  =1$ and $cm\equiv cn \pmod{3}$. Suppose (2) fails. The reader can check that $3\nmid c,$ in which case $m\equiv n \pmod{3}$. But $x/c\equiv y/c \pmod{3}$ and the primitivity of $(x,y)$ imply $c=1$.  
\end{proof}

\begin{example}
\label{ex}Using Theorem \ref{PrimeCheck}, the reader can check that the
following orbits of period $2n$ are primitive:

$\bullet$ $\ n=2k+1,k\geq1:(k-1,k+2)$

$\bullet$ $\ n=2:(1,1)$

$\bullet$ $\ n=4k+4,k\geq1:(2k-1,2k+5)$

$\bullet$ $\ n=4k+10,k\geq1:(2k-1,2k+11)$.
\end{example}

\noindent Since $\mathcal{P}(n)$ tells us there are \emph{no}
primitive orbits of period $2,$ $8,$ $12$ or $20,$ Example \ref{ex} exhibits
a primitive orbit of every possible even period.

\section{CONCLUDING REMARKS}

Many interesting open questions remain; we mention three:

\smallskip

\noindent(1) What can be said if the equivalence relation on the set of all
periodic orbits defined above is defined more restrictively? For example,
one could consider an equivalence relation in which equivalent orbits have
cycles of incidence angles that differ by a \emph{cyclic }permutation.

\smallskip

\noindent(2) Every isosceles triangle admits a period 4 orbit resembling
(1,1) and every acute triangle admits an orbit of period 6 resembling (0,3).
Empirical evidence suggests that every acute isosceles triangle with base
angle at least 54 degrees admits an orbit of period 10 resembling (1,4).
Thus we ask: To what extent do the results above generalize to acute
isosceles triangles?

\smallskip

\noindent(3) Arbitrarily label the sides of the triangle $0,1,2$ and consider the sequence of integers modulo 3 given by the successive bounces of a billiards trajectory.  Clearly periodic trajectories yield periodic sequences.  For example, the sequence $01020102\ldots$ is given by the period 4 orbit (1,1).  If $\{a_n\}$ is a periodic mod 3 sequence, is $\{a_n\}$ given by some billiards trajectory?

\section{ACKNOWLEDGEMENTS}

This project emerged from an undergraduate research seminar directed by
Zhoude Shao and the second author during the spring of 2003. Student
participants included John Gemmer, Sean Laverty, Ryan Shenck, Stephen Weaver,
and the first author. To assist us computationally, Stephen Weaver created
his \textquotedblleft Orbit Tracer\textquotedblright\ software \cite{W1},
which generated copious experimental data and produced the diagrams in
Figure 10 above. Dennis DeTurck suggested we consider the general billiards
problem and consulted with us on several occasions.  Numerous persons read the manuscript and offered helpful suggestions at various stages of its development. These include Annalisa Crannell, Doris Schattschneider, Jim Stasheff, Doron Zeilberger, and the referees.  We thank each of these individuals for their contributions.

Andrew Baxter received his B.A. in Mathematics from Millersville University in 2005.  This article, his first published work, is based on his undergraduate research thesis which won the MAA Eastern Pennsylvania and Delaware chapter Student Paper Competition.  He is currently pursuing his Ph.D. in combinatorics at Rutgers University under the supervision of Doron Zeilberger.  His wife Kristen is also pursuing her Ph.D. at Rutgers in Classical Studies.

Ron Umble is a professor of mathematics at Millersville
University of Pennsylvania, where he has been
a faculty member since 1984. His research interests
are in algebraic topology and the theory of generalized
operads, in particular. He has directed numerous
undergraduate research projects and has coauthered
three published papers with students (including this
one). He is a member of the Pennsylvania Zeta
chapter (Temple University) of Pi Mu Epsilon and
was the 1972 winner of the Zeta chapter's student
paper competition.

\end{document}